\documentclass[11pt,letterpaper]{article}

\usepackage[utf8]{inputenc}
\usepackage[T1]{fontenc}
\usepackage{lmodern}
\usepackage{fullpage}

\usepackage{amsthm,amsmath,amssymb}
\usepackage{amsfonts}

\usepackage{graphicx}
\usepackage{floatrow}

\usepackage{hyperref}
\usepackage{amsrefs}   

\def\Z{\mathbb{Z}}

\def\GG{\mathcal{G}}

\makeatletter
\def\proof{\@ifnextchar[
  {\@xproof}{\@proof}}
\def\endproof{\unskip\kern 10\p@ \begingroup \unitlength\p@
  \linethickness{.4\p@} \framebox(6,6){} \endgroup \endtrivlist }
\def\@proof{\trivlist \item[\hskip\labelsep {\em Proof.}  ]
  \ignorespaces }
\def\@xproof[#1]{\trivlist \item[\hskip\labelsep{\em #1.}]
  \ignorespaces }
\makeatother

\theoremstyle{plain}
\newtheorem{theorem}{Theorem}

\newtheorem{proposition}[theorem]{Proposition}
\newtheorem{lemma}[theorem]{Lemma}

\newtheorem*{observation-sn}{Observation}

\theoremstyle{plain} 
\newcommand{\thistheoremname}{}
\newtheorem{genericthm}[theorem]{\thistheoremname}

\theoremstyle{remark}

\theoremstyle{plain}

\def\ww{{\mathcal W}}

\usepackage{mathrsfs}

\def\ww{{\mathcal W}}
\def\ignore#1{{}}

\def\gen{{\text{\rm gen}}}

\title{\bf Embeddability of arrangements of pseudocircles \\ and graphs on surfaces}

\author{\'Eric Colin de Verdi\`ere \thanks{Universit\'e Paris-Est, LIGM, CNRS, ENPC, ESIEE Paris, UPEM, Marne-la-Vall\'ee, France. {\em E-mail:} {\tt eric.colindeverdiere@u-pem.fr}.}
\and
Carolina Medina\thanks{Instituto de Física, UASLP. 
San Luis Potosí, Mexico, 78000.
{\em E-mail:} {\tt cmedina@math.ucdavis.edu}.}
\and
Edgardo Roldán-Pensado\thanks{Centro de Ciencias Matemáticas, UNAM campus Morelia. Morelia, Michoacán, Mexico, 58190. {\em E-mail:} {\tt e.roldan@im.unam.mx.}}
\and
Gelasio Salazar\thanks{ Instituto de Física, UASLP. 
San Luis Potosí, Mexico, 78000. {\em E-mail:} {\tt gsalazar@ifisica.uaslp.mx}. 
}}

\begin{document}


\maketitle

\begin{abstract}
  A pseudocircle is a simple closed curve on some surface; an arrangement of pseudocircles is a collection of pseudocircles that pairwise intersect in exactly two points, at which they cross.  Ortner proved that an arrangement of pseudocircles is embeddable into the sphere if and only if all of its subarrangements of size at most four are embeddable into the sphere, and asked if an analogous result holds for embeddability into orientable surfaces of higher genus. We answer this question positively: An arrangement of pseudocircles is embeddable into an orientable surface of genus~$g$ if and only if all of its subarrangements of size at most $4g+4$ are.  Moreover, this bound is tight. We actually have similar results for a much general notion of arrangement, which we call an \emph{arrangement of graphs}.
%
\end{abstract} 


\section{Introduction}\label{sec:intro}

The starting point of this work is motivated by Ortner~\cite{ortner}: He proved that an arrangement of pseudocircles is embeddable into the sphere if and only if all of its subarrangements of size at most four are embeddable into the sphere, and asked if an analogous result held for embeddability into surfaces of higher genus.  We answer this question positively, and in fact prove a similar result for the more general notion of \emph{arrangement of graphs}.

We briefly recall some standard notions of topological graph theory; see Mohar and Thomassen~\cite{moharthomassen} for details.  All surfaces under consideration are orientable; we denote by $S_g$ the (orientable) surface of genus~$g$. If $G$ is a graph and $v$ is a vertex of $G$, then a {\em rotation} of $v$ is a cyclic ordering of the edges incident with $v$. A {\em rotation system} of $G$ is a collection of rotations of all vertices of $G$. 

\medskip
\noindent{\bf Remark.} {\em Throughout this work, every graph under consideration is implicitly assumed to be equipped with a rotation system. {(This enhanced notion of a graph is also called a {\em combinatorial map} in the literature). If} $G$ is a graph and $H$ is a subgraph of $G$, then it is implicitly assumed that the rotation system of $H$ is the one naturally inherited from the rotation system of $G$.}

\medskip
\medskip

We define an {\em arrangement of graphs}  $\GG$ to be a collection $(G_0, \ldots, G_n)$ of connected subgraphs of a  graph $\bigcup_{i=0}^n G_i$ such that $G_0$ has at least one vertex in common with each of $G_1,\ldots,G_n$.  For brevity, we use $\bigcup\GG$ to denote $\bigcup_{i=0}^n G_i$. We emphasize that $\bigcup\GG$ (as all graphs under consideration) is endowed with a rotation system, from which $G_i$ inherits a rotation system, for $i=0,\ldots,n$. The {\em size} of $\GG$ is $n+1$, the number of graphs in the arrangement.  A {\em subarrangement} of $\GG$ is a subcollection of $\GG$ that contains $G_0$. Thus every subarrangement of $\GG$ is an arrangement.

An {\em embedding} of a graph $G$ on an orientable surface~$S$ is a drawing of $G$ on $S$ without crossings, such that the clockwise ordering of the edges around each vertex $v$ in this drawing matches the rotation of $v$ in the rotation system of $G$. An {\em embedded graph} is a graph with a given embedding on a surface. A graph $G$ is {\em embeddable} into a surface $S$ if there is an embedding of $G$ on $S$. 

Our main result is the following.
\begin{theorem}\label{thm:gen}
  Let $g\ge0$.  An arrangement of graphs is embeddable into~$S_g$ if and only if all its subarrangements of size at most $4g+5$ are embeddable into $S_g$.  Moreover, the bound of $4g+5$ is tight.
\end{theorem}
In Theorem~\ref{thm:gen}, the requirement that some graph~$G_0$ intersects all other graphs in the arrangement may seem superfluous.  However, as we argue in Section~\ref{sec:concludingremarks}, some condition along these lines is required in order to obtain a result in this spirit.

Let us recast the notion of arrangement of pseudocircles in our terminology.  A {\em pseudocircle} is a simple closed curve in a surface.  Gr\"unbaum~\cite{grunbaum} defined an \emph{arrangement of pseudocircles} (actually, he used the terminology \emph{arrangement of curves}) to be a family of pseudocircles embedded on a surface such that any two pseudocircles intersect in exactly two points, at which they cross.   In our language, an arrangement of pseudocircles is thus an arrangement of graphs $\GG=(G_0, \ldots, G_n)$ such that each $G_i$ is a cycle and, for each $i\neq j$, $G_i$ and~$G_j$ intersect at exactly two vertices of~$\bigcup\GG$, where they cross.  The main result by Ortner~\cite{ortner} is the following.

\begin{theorem}[Ortner~\cite{ortner}*{Theorem~10}]\label{thm:ortner}
  An arrangement of pseudocircles is embeddable into~$S_0$ if and only if all its subarrangements of size at most four are embeddable into~$S_0$.
\end{theorem}

At the end of his article, Ortner asks whether such a result can be generalized to arbitrary surfaces.  Theorem~\ref{thm:gen} already answers this question positively.  We actually prove the following result, with a sharp bound on the size of the subarrangements:
\begin{theorem}\label{thm:ortner-gen}
  Let $g\ge0$.  An arrangement of pseudocircles is embeddable into~$S_g$ if and only if all its subarrangements of size at most $4g+4$ are embeddable into $S_g$.  Moreover, the bound of $4g+4$ is tight.
\end{theorem}

Depending on the context, the notion of arrangement of pseudocircles can be relaxed in several ways; see for instance~\cite{ortner2,kang,anppss,linortner3}.  Arrangements of other objects, such as line segments~\cite{egppss} (or Jordan arcs~\cite{ede1}), can also be naturally regarded as special cases of arrangements of graphs, and so our Theorem~\ref{thm:gen} also applies to them.


To explain our proof strategy, let us first define the {\em genus} $\gen(G)$ of a graph $G$ to be the smallest integer $g$ such that $G$ embeds in~$S_g$.  Such an embedding of~$G$ is necessarily cellular, and one can obtain it by attaching a disk to each boundary walk of~$G$ (those boundary walks are well defined by the rotation system). If $G$ has vertex set $V$, edge set $E$, and set $\ww$ of boundary walks, then, by Euler's formula, we have the following~\cite{moharthomassen}*{Section 4.1}:
\begin{equation}
\gen(G)=(1/2)(2-|V |+|E|-|\ww|).\label{eq:euler}
\end{equation}

In Section~\ref{sec:proofmain}, we shall prove the following proposition.
\begin{proposition}\label{P:main}
Let $\GG$ be an arrangement of graphs, and let $h:=\gen(\bigcup\GG)$. Then for each $g\in \{0,\ldots,h\}$ there is a subarrangement $\GG_g$ of $\GG$, of size at most $4g+1$, such that $\gen(\bigcup\GG_g) \ge g$. Moreover, if $\GG$ is an arrangement of pseudocircles, then this bound of $4g+1$ can be improved to $4g$ for every $g\in\{1,\ldots,h\}$.
\end{proposition}

In Section~\ref{sec:proofortnerg}, we prove our Theorems~\ref{thm:gen} and~\ref{thm:ortner-gen}.  As we will see, the upper bounds of $4g+5$ and $4g+4$ follow as an easy corollary of the above proposition; we use specific constructions to prove their tightness.  As we will see, the tightness of the bound $4g+5$ in Theorem~\ref{thm:gen} is witnessed by arrangements of graphs which are quite specific: Each graph is a cycle, and two distinct graphs are either disjoint or have exactly two intersection points, at which they cross.

\bigskip
\noindent{{\bf Remark.} \emph{In general, when considering graphs, one allows loops and multiple edges. However, it is not hard to see that proving Theorem~\ref{thm:gen} reduces to proving it in the case of graphs without loops (simply by subdividing each loop once).  So, henceforth, we assume that all graphs under consideration are without loops.}}

\section{Proof of Proposition~\ref{P:main}}\label{sec:proofmain}

A {\em face} of a graph $G$ embedded on a surface $S$ is a connected component of $S\setminus G$. We remark that a face $f$ is an open subset of $S$, that is, the vertices and edges that lie on boundary of $f$ are not part of $f$. {Every face is bounded by one or several walks of the graph, called \emph{boundary walks}, which are determined by the rotation system of $G$.}

We follow the convention that as we traverse a boundary walk of an embedding of a graph $G$, the face lies at our right-hand side. Thus, for instance, if $G$ is just a cycle $v_0 e_1 v_1 \cdots e_m v_0$ in the sphere, then $G$ has two faces, one of which has boundary walk $v_0 e_1 \cdots e_m v_0$, and the other has boundary walk $v_0 e_m \cdots e_1 v_0$. We note that this convention would be ambiguous in the presence of loop edges, but by the remark at the end of the previous section this is not an issue in our context since all graphs under consideration are loopless.

Let $G$ be a graph, and let $J$ be a subgraph of $G$. If $G$ is embedded on a surface $S$, then $J$ naturally inherits an embedding on $S$ from the embedding of $G$, and so we will implicitly regard $J$ also as an embedded graph.

Let $G$ be an embedded graph, and let $J$ be a subgraph of $G$.  
Let $f$ be a face of $J$, and let $e$ be an edge of $G\setminus J$ contained in $f$, incident with a vertex $v$ in a boundary walk $W$ of $f$. If as we traverse $W$ we encounter $e$ leaving $v$ at our right-hand side, then we say that $e$ {\em attaches to $W$ at $v$}. Note that $e$ may attach to (at most) two distinct walks, and if it attaches to two walks, then it may do so at the same vertex or at different vertices. 

An instance in which an edge attaches to two distinct walks at the same vertex is given by the following example. Let $G$ be a cellularly embedded graph in the torus with only one vertex $v$ and two loop-edges $e$ and $e'$, and let $J=G-e$. Then $J$ has two boundary walks, and as we traverse either of these boundary walks we encounter $e$ leaving $v$ at our right-hand side. 

Continuing with the discussion, let $G$ be an embedded graph, let $J$ be a subgraph of $G$, and let $f$ be a face of $J$. Let $W,W'$ be distinct boundary walks of $f$. A walk $U=v_0e_1v_1\ldots e_m v_m$ of $G$ is a $WW'$-{\em walk} if $e_1,v_1,\ldots,v_{m-1}, e_m$ are in the interior of~$f$, $e_1$ attaches to $W$ at $v_0$, and $e_m$ attaches to $W'$ at $v_m$. Thus a $WW'$-walk has its (not necessarily distinct) endvertices in $J$, and is otherwise contained in $G\setminus J$. If $U$ is either a path or a cycle, we say that it is a $WW'$-{\em path} {\em contained in $f$}.

The proof of Proposition~\ref{P:main} relies crucially on the following two lemmas. 

\begin{lemma}\label{L:2bd}
Let $\GG=(G_0,\ldots,G_n)$ be an arrangement of graphs such that $\bigcup\GG$ is cellularly embedded on a surface~$S$, and let $H$ be a
  subgraph of~$\bigcup\GG$ that contains the vertex set of~$G_0$.  Assume that some face~$f$ of~$H$ has at least two
  boundary walks.  Then for some distinct boundary walks $W,W'$ of $f$ there is a $WW'$-path contained in $f$, whose edges are contained in the union of two graphs in~$\GG$. 
\end{lemma}

\begin{lemma}\label{L:1bd}
  Let $\GG=(G_0,\ldots,G_n)$ be an arrangement of graphs such that $\bigcup\GG$ is cellularly embedded on a surface~$S$, and let $H$ be a subgraph of~$\bigcup\GG$ that contains the vertex set of~$G_0$.  Assume that some face~$f$ of~$H$ has positive genus and a single boundary walk. Then there is a subgraph $Q$ of\, $\bigcup\GG$ that lies in $f\cup \partial f$, whose edges belong to the union of two graphs in~$\GG$, and such that $Q$ is either:
\begin{enumerate}
\item a path with endpoints in~$\partial f$ that is otherwise disjoint from $\partial f$, and does not separate $f$; or
\item a cycle with one vertex in $\partial f$ that is otherwise disjoint from $\partial f$, and does not separate $f$; or
\item the union of a cycle $C$ contained in the interior of~$f$, that does not separate $f$, and of a path that has one endpoint in $C$ and one endpoint in $\partial f$, and is otherwise disjoint from $C$ and $\partial f$. 
\end{enumerate}
\end{lemma}

\begin{proof}[Proof of Proposition~\ref{P:main}, assuming Lemmas~\ref{L:2bd} and~\ref{L:1bd}]
Let us consider an embedding of~$\bigcup\GG$ on~$S_h$; since $h=\gen(\bigcup\GG)$, the embedding is necessarily cellular.  We prove by induction on $g=0,\ldots,h$ that there is an embedded subgraph $H_g$ of $\bigcup\GG$ such that:
\begin{itemize}
\item $H_g$ contains all the vertices of $G_0$;
\item $H_g$ has genus~$g$;
\item $H_g$ has a single face in~$S_h$, with a single boundary walk;
\item $H_g$ is contained in a subarrangement of~$\GG$ of size at most $4g+1$.
\end{itemize}
This will show the first part of the proposition: Indeed, let $\GG_g$ be a subarrangement of~$\GG$ of size at most $4g+1$ containing~$H_g$; we have $g=\gen(H_g) \le \gen(\bigcup\GG_g)$.

For the base case, it suffices to take $H_0$ to be a spanning tree of~$G_0$.  For the inductive step, let $0<g\le h$ and suppose the statement holds for $g-1$. Because $g-1<h$, the unique face~$f$ of $H_{g-1}$ has a single boundary walk. Moreover, $f$ has positive genus: since $S_h$ has genus $h$ and $H_{g-1}$ has genus $g-1$, it follows that the genus of $f$ is $h-(g-1)$.

Let us first apply Lemma~\ref{L:1bd} to~$H:=H_{g-1}$ and to its unique face~$f$; we obtain a subgraph~$Q$ of~$\bigcup\GG$ lying in $f\cup\partial f$, whose edges belong to the union of two graphs in~$\GG$, having one of the three specific structures mentioned in Lemma~\ref{L:1bd}.  Regardless of that structure, the graph $H'_{g-1}:=H_{g-1}\cup Q$ has a single face, with exactly two boundary walks.  Moreover, if $H'_{g-1}$ was obtained from~$H_{g-1}$ by adding $k$~edges, then $k-1$ vertices were added.  Thus, by Equation~\ref{eq:euler}, the genus of~$H'_{g-1}$ equals that of~$H_{g-1}$, namely~$g-1$.

Let us now apply Lemma~\ref{L:2bd} to~$H:=H'_{g-1}$ and to its unique face~$f$; if $W$ and~$W'$ are the two boundary walks of~$f$, we obtain a $WW'$-path~$P$ whose edges are contained in the union of two graphs in $\GG$. We finally let $H_g:=H_{g-1}'\cup P$.  This graph has a single face, with exactly one boundary walk.  Moreover, if $H_g$ was obtained from~$H'_{g-1}$ by adding $k$~edges, then $k-1$ vertices were added.  Thus, by Equation~\ref{eq:euler}, the genus of~$H_g$ equals that of~$H'_{g-1}$ plus one, namely~$g$.

Finally, recall that $H_{g-1}$ is contained in a subarrangement of~$\GG$ of size at most $4g-3$.  Thus, by construction, $H_g$ is contained in a subarrangement of~$\GG$ of size at most $4g+1$.  The proof of the induction step is complete.

We proceed analogously if $\GG$ is induced from an arrangement of pseudocircles, but in this case we proceed by induction on $g=1,\ldots,h$.  The induction hypothesis is the same, with $4g+1$ replaced by $4g$; the proof of the inductive step is identical.  Let us prove the base case $g=1$.  Theorem~\ref{thm:ortner} implies that (since $\bigcup\GG$ is not embeddable into the sphere, as $\gen(\bigcup\GG)\ge1$) $\GG$ has a subarrangement $\GG_1$ of size at most four, such that $\bigcup\GG_1$ is not embeddable into the sphere, that is, $\gen(\bigcup\GG_1)\ge 1$.  Applying to~$\GG_1$ the same construction as the induction step above, we obtain a subgraph~$H_1$ of~$\bigcup\GG_1$ of genus exactly one, embedded on~$S_h$ with a single face, having a single boundary walk.  Of course, since~$\GG_1$ has size four, the graph~$H_1$ is contained in a subarrangement of~$\GG$ of size at most four, as desired.
\end{proof}


\begin{proof}[Proof of Lemma~\ref{L:2bd}]
Let $f$ be a face of $H$ whose set $\ww=\{W_1,\ldots,W_r\}$ of boundary walks has size at least $2$. The cellularity of $\bigcup\GG$ implies that for some two distinct $W_j,W_k$ in $\ww$ there exists a $W_jW_k$-path contained in $f$. If there is such a $W_jW_k$-path $U$ with only one or two edges then we are done, as the edges of $U$ are obviously contained in the union of two graphs of $\GG$. Thus we assume that every $W_jW_k$-path has at least three edges. In particular, there are vertices of $\bigcup\GG$ contained in $f$.

We say that a vertex $v$ contained in $f$ has {\em colour} $\ell\in\{1,\ldots,r\}$ if there is a path $P$ with the following properties: (i) $P$ starts at $v$, its final edge attaches to $W_\ell$, and except for this attachment, $P$ is contained in $f$; and (ii) there is a $G_i\in\GG$ that contains all the edges of $P$. 

We note that each vertex contained in $f$ has at least one colour. Indeed, let $w$ be a vertex contained in $f$. Since $H$ contains all the vertices of $G_0$, it follows that there is an $i\in\{1,\ldots,n\}$ such that $w$ is in $G_i$. Since $G_i$ is connected and it contains at least one vertex of $G_0$, it follows that there is a path contained in $G_i$ that has $w$ as an endpoint and whose other endpoint is in $G_0$. A shortest path $P$ with this property has the endpoint in $G_0$ necessarily in the boundary of $f$, and so the final edge of $P$ attaches to $W_\ell$ for some $\ell\in \{1,\ldots,r\}$.

We also note that ($*$) if $e=uv$ is an edge such that $e,u$, and $v$ are all contained in $f$, then $u$ and $v$ have at least one common colour. To see this, first we note that since $e$ and its endvertices are contained in $f$, then there is an $i\in \{1,\ldots,n\}$ such that $e$ (and hence also $u$ and $v$) belongs to $G_i$. Using the same arguments as in the previous paragraph, it follows that there is a path $Q$ contained in $G_i$, whose first edge is $e$ (the startpoint of $Q$ is either $u$ or $v$) and that attaches to $W_\ell$ for some $\ell\in\{1,\ldots,r\}$. Therefore $u$ and $v$ have the common colour $\ell$.

We shall show that some vertex contained in $f$ has at least two distinct colours. Note that this completes the proof: if $j\neq k$ are both colours of $v$, then it follows that there is a $W_jW_k$-walk, and hence a $W_jW_k$-path, contained in the union of two graphs in $\GG$.

We recall from the first paragraph of this proof that there is a $W_jW_k$-path $U=v_0e_1 v_1\cdots e_mv_m$ contained in $f$, for some distinct $W_j,W_k\in\ww$, where $m\ge 3$. By ($*$), for each $i=1,\ldots,m-2$ the vertices $v_i$ and $v_{i+1}$ have a colour in common. Since $v_1$ has color $j$ and $v_{m-1}$ has colour $k$, it follows that at least one vertex in $\{v_1,\ldots,v_{m-1}\}$ has at least two colours.
\end{proof}

The heart of the proof of Lemma~\ref{L:1bd} is the following lemma, which is reminiscent of Thomassen's \emph{3-path condition}~\cite{thomassen}; see also Mohar and Thomassen~\cite{moharthomassen}*{Chapter~4}.

\begin{lemma}\label{lem:nonsep}
Let $\GG=(G_0,\ldots,G_n)$ be an arrangement of graphs such that $\bigcup\GG$ is cellularly embedded on a surface $S$ of positive genus, and $G_0$ consists of a single vertex $v$. Then there is a non-separating cycle in $\bigcup\GG$ contained in the union of two graphs in $\GG$.
\end{lemma}

\begin{proof}
  We use cellular homology over~$\Z/2\Z$, but we make the proof self-contained.  Let $G=(V,E)$ be the graph equal to~$\bigcup\GG$.  A subset $E'$ of~$E$ is a \emph{boundary} if the faces of~$G$ can be colored in black and white in such a way that $E'$ is exactly the set of edges of~$G$ incident with one black face and one white face.  

{We remark that, in contrast to homology theory, here a \emph{cycle} is a closed walk without repeated vertices.}

We have the following (standard) properties:
  \begin{enumerate}
  \item The symmetric difference of two boundaries is a boundary.
  \item Let $E'$ be the edge set of a cycle~$C$ in~$G$.  Then $E'$ is a boundary if and only if $C$ is separating.
  \end{enumerate}
  
  Let $v_0e_1v_1\ldots e_m v_m$ (with $v_m=v_0$) be a non-separating cycle~$C$ in~$\bigcup\GG$.  For each $i=0,\ldots,m-1$, let us choose a path~$P_i$ from $v$ to~$v_i$ in some graph~$G_p$ of the arrangement, such that $G_p$ also contains~$e_i$.  Finally, let $L_i$ be the set of edges of~$E$ that appear in the closed walk $P_i\cdot e_{i+1}\cdot\overline{P_{i+1}}$ (where $\overline{P_{i+1}}$ denotes the reversal of~$P_{i+1}$) an odd number of times.  By construction, for each~$i$ there exist two graphs in~$\GG$ whose union contains all the edges in~$L_i$.

We remark that the symmetric difference of the $L_i$ equals precisely the edge set of~$C$, which is not a boundary since $C$ is non-separating (Property~2). By Property~1, there is an $i\in\{0,\ldots,m-1\}$ such that $L_i$ is not a boundary.

{We now claim} that each vertex in the graph $(V,L_i)$ has even degree. To see this, let $u$ be a vertex in $(V,L_i)$ with degree greater than $0$. Thus $u$ appears in the walk $W_i:=P_i\cdot e_{i+1}\cdot\overline{P_{i+1}}$. Define the {\em weight} of each edge $e$ in $W_i$ incident with $u$ as the number of times that $e$ appears in $W_i$. Since $W_i$ is closed, the sum of the weights of all edges in $W_i$ incident with $u$ is even, and so the number of edges in $W_i$ incident with $u$ that have odd weight must be even. Since the number of edges in $W_i$ incident with $u$ that have odd weight is precisely the degree of $u$ in $(V,L_i)$, this proves the claim.

{It follows} that $L_i$ is the disjoint union (and thus, the symmetric difference) of edge sets of cycles; for the same reason as above, one of these edge sets of cycles is not a boundary.  By Property~2 again, the corresponding cycle is non-separating; like~$L_i$, it is contained in the union of some two graphs in~$\GG$.
\end{proof}

\begin{proof}[Proof of Lemma~\ref{L:1bd}]
  Let $S'$ be the surface that results by contracting $S\setminus f$ to a single point. Then $\GG$ naturally induces an arrangement $(G_0',G_1',\ldots,G_n')$ in $S'$, where $G_0'$ consists of a single vertex $v$. Note that the graphs in $\GG$ that do not have any edge in $f$ also get collapsed to a single vertex. Also note that the cellularity of $\GG$ implies that $\GG'$ is cellularly embedded in $S'$. By Lemma~\ref{lem:nonsep}, there exist (not necessarily distinct) $G_i',G_j'$ such that $G_i'\cup G_j'$ contains a non-separating cycle $C$. 

It is easy to see that if $C$ contains $v$, then one of the first two outcomes in Lemma~\ref{L:1bd} must hold. Suppose finally that $C$ does not contain $v$. Every graph in the arrangement contains $v$ (as $G_0'$ consists only of $v$), and in particular $G_i'$ contains $v$. This implies that there is a path $P$ from $C$ to $v$ that is contained in $G_i'$, such that $P$ only intersects $C$ at the initial vertex of $P$. In this case we have the third outcome in Lemma~\ref{L:1bd}. \end{proof}

\section{Proof of Theorems~\ref{thm:gen} and~\ref{thm:ortner-gen}}\label{sec:proofortnerg}

\subsection{Upper bounds}

First we prove the upper bound of~$4g+5$, for arbitrary arrangements of graphs (Theorem~\ref{thm:gen}). We only need to show the ``if'' part in the theorem, as the ``only if'' part is trivial. We prove the contrapositive statement: we let $\GG$ be an arrangement of graphs not embeddable into~$S_g$, and show that $\GG$ has a subarrangement of size at most $4g+5$ that is not embeddable into $S_g$.

Let $h:=\gen(\bigcup\GG)$; we have $h\ge g+1$ because $\GG$ is not embeddable into~$S_g$.  By Proposition~\ref{P:main}, $\GG$ has a subarrangement $\GG_{g+1}$ of size at most $4(g+1)+1=4g+5$, such that $\gen(\bigcup\GG_{g+1}) \ge g+1$. This implies that the subarrangement~$\GG_{g+1}$ of~$\GG$, which has size at most~$4g+5$, is not embeddable into $S_g$, as desired.

To prove the upper bound of~$4g+4$ in the case of arrangements of pseudocircles (Theorem~\ref{thm:ortner-gen}), we follow the same arguments, using the stronger bound $4g$ guaranteed by Proposition~\ref{P:main}.

\subsection{Tightness of the bounds}

\begin{center}
\begin{figure}[ht!]
\includegraphics[width=14cm]{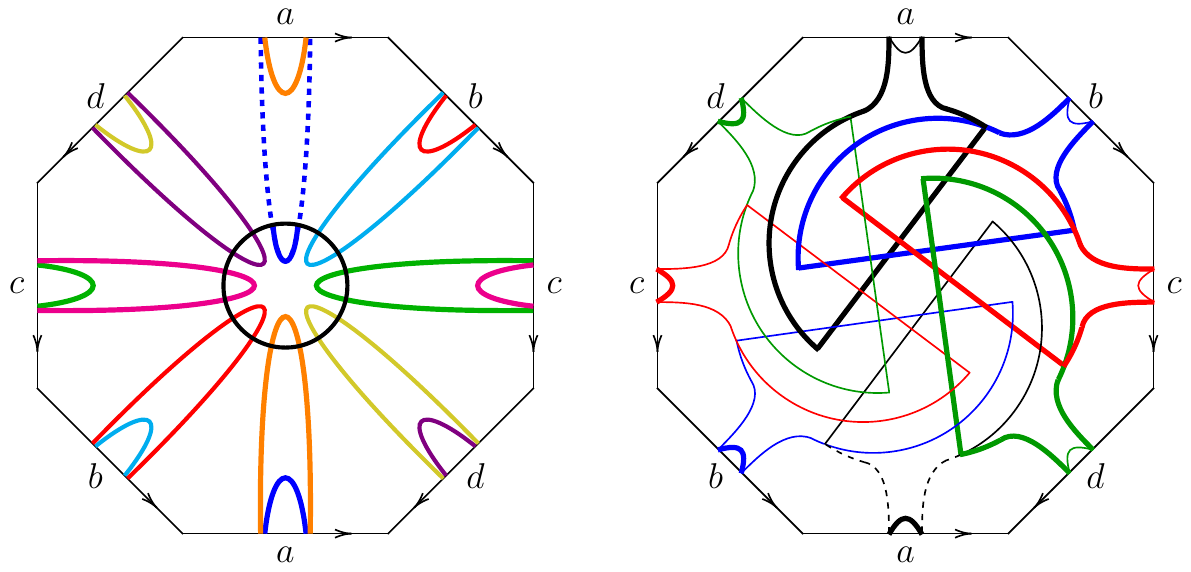}
\caption{Illustration of the tightness of the bounds in Theorems~\ref{thm:gen} (left) and~\ref{thm:ortner-gen} (right).}
\label{fig:figex1}
\end{figure}
\end{center}

We first prove that the upper bound of~$4g+5$, for arbitrary arrangements of graphs, is tight (Theorem~\ref{thm:gen}).  For this, we refer the reader to Figure~\ref{fig:figex1}. On the left hand side we have an arrangement of graphs~$\GG$ of size $9$ in the double torus $S_2$; each graph is a cycle, and two distinct graphs are either disjoint, or have exactly two intersection points, at which they cross. Since $\bigcup\GG$ is cellularly embedded in $S_2$, it follows that $\bigcup\GG$ is not embeddable into the torus. We claim that every subarrangement of $\GG$ of size $8$ is embeddable into the torus.

To see this first we note that every subarrangement of $\GG$ must contain the graph $G_0$ at the center of the polygon, as this is the only graph that intersects all the other graphs. Thus a subarrangement of size $8$ is obtained by removing any of the other graphs. We claim that not only every subarrangement $\GG'$ of $\GG$ of size $8$ is embeddable into the torus, but that it suffices to remove two particular edges of an arbitrary graph (distinct from $G_0$) in order to obtain an embedded graph with genus $1$. 

Indeed, let $G$ be the subgraph of $\bigcup\GG$ obtained by removing the two dashed edges on the left hand side of Figure~\ref{fig:figex1}. Thus $G$ has the same number of vertices as $\bigcup\GG$, and two fewer edges than $\bigcup\GG$. It is easy to check that $G$ has the same number of boundary walks as $\bigcup\GG$. Using \eqref{eq:euler}, it follows that $\gen(G)=\gen(\bigcup\GG)-1=1$, as claimed. 

We finally note that this construction is easily extended to get, for each $g\ge 0$, an arrangement of graphs (all of them being cycles) of size $4g+5$ that is not embeddable into the surface $S_g$ of genus $g$, and all of whose subarrangements of size $4g+4$ are embeddable into $S_g$.

\bigskip

The tightness of the bound of~$4g+4$ for arrangements of pseudocircles (Theorem~\ref{thm:ortner-gen}) is illustrated on the right hand side of Figure~\ref{fig:figex1}. This figure shows an arrangement of pseudocircles $\GG$ of size $8$ in $S_2$. Since $\bigcup\GG$ is cellularly embedded in $S_2$, it follows that $\bigcup\GG$ is not embeddable into the torus. Now remove from $\bigcup\GG$ the two dashed edges shown in the figure. It is easy to verify that the resulting embedded graph $G$ has two fewer edges and the same number of vertices and boundary walks as $\bigcup\GG$. Using~\ref{eq:euler}, it follows that $\gen(G)=\gen(\bigcup\GG)-1=1$. The two deleted edges belong to the same pseudocircle, and so it follows that the subarrangement of $\GG$ obtained by removing this pseudocircle is embeddable into the torus. By the symmetry of the construction, every subarrangement of $\GG$ of size $7$ is embeddable into the torus.

This construction is easily extended to give, for each $g \ge 0$, an arrangement of pseudocircles of size $4g+4$ that is not embeddable into the surface $S_g$ of genus $g$, and all of whose subarrangements of size $4g+3$ are embeddable into $S_g$.

\section{Concluding remarks}\label{sec:concludingremarks}

To obtain a result along the lines of Theorems~\ref{thm:gen} or~\ref{thm:ortner-gen}, it is natural to ask if it is absolutely necessary to require that there is a graph in the arrangement intersecting all other graphs in the collection. To answer this question, we note that it is necessary to require some sort of condition along these lines. Indeed, as observed by Ortner~\cite{ortner}*{Figure 16}, there exist arbitrarily large collections of pseudocircles (whose union is connected) that cannot be embedded into a sphere, and yet the removal of any pseudocircle leaves an arrangement that can be embedded into a sphere.

On the other hand, in order to have some version of Theorem~\ref{thm:gen}, it is not strictly necessary to have a single graph intersecting all the others; our techniques and arguments are readily adapted under the assumption that there is a subcollection of bounded size that gets intersected by all other graphs. More precisely, let us define an $m$-{\em arrangement} of graphs as a collection in which there is a subcollection~$B$ of size (at most) $m$ such that every graph intersects at least one graph in~$B$, and the union of the graphs is connected. It is easy to verify that there is a choice of at most $2(m-1)$ graphs whose union with~$B$ is connected; this union is, of course, intersected by every other graph.   Thus the graphs in an $m$-arrangement naturally induce an arrangement of graphs $(G_0,\ldots,G_n)$, where $G_0$ is the union of at most $3m-2$ graphs. The following is then an easy consequence of Proposition~\ref{P:main}.

\begin{theorem}
  An $m$-arrangement of pseudocircles is embeddable into $S_g$ if and only if all of its subarrangements of size at most $4g+3m+2$ are embeddable into $S_g$.
\end{theorem}

\section*{Acknowledgements}

We thank two anonymous referees for many helpful suggestions and corrections to an earlier version of this paper. The first author was supported by grant ANR-17-CE40-0033 of the French National Research Agency ANR (SoS project). The second author is currently supported by a Fulbright Visiting Scholar grant at UC Davis. The second and fourth authors were supported by Conacyt under Grant 222667, and by FRC-UASLP. The third author was supported by Conacyt under Grant 166306 and by PAPIIT-UNAM IA102118.


\end{document}